\documentclass[12pt,leqno]{amsart}%
\usepackage{amsmath}
\usepackage[margin=1.5in]{geometry}
\usepackage{graphicx}
\usepackage{amsfonts}
\usepackage{amssymb}%
\usepackage[utf8]{inputenc}
\usepackage{comment}
\usepackage{hyperref}
\usepackage{mathtools}
\usepackage{ dsfont }

\newtheorem{theorem}{Theorem}[section]
\theoremstyle{plain}

\newtheorem{claim}{Claim}

\newtheorem{lemma}{Lemma}[section]

\newtheorem{proposition}{Proposition}[section]

\numberwithin{equation}{section}
\theoremstyle{definition}

\theoremstyle{remark}
\newtheorem{remark}{Remark}[section]

\begin{document}
\title[Lagrangian mean curvature equation]{A note on the two dimensional Lagrangian mean curvature equation}
\author{Arunima Bhattacharya}
\address{Department of Mathematics\\
University of Washington, Seattle, WA 98195, U.S.A.}
\email{arunimab@uw.edu}

\begin{abstract}
 In this note, we use Warren-Yuan's \cite{WY2d} super isoperimetric inequality on the
level sets of subharmonic functions, which is available only in two dimensions, to 
 derive a modified Hessian bound for solutions of the two dimensional Lagrangian mean curvature equation. We assume the Lagrangian phase to be supercritical with bounded second derivatives. Unlike in \cite{AB}, the simplified approach in this paper does not require the Michael-Simon mean value and Sobolev inequalities on generalized submanifolds of $\mathbb{R}^n$ \cite{MS}. 

\end{abstract}

\maketitle

\section{ Introduction}

We study the two dimensional Lagrangian mean curvature equation
\begin{equation}
    \arctan \lambda_{1}+\arctan \lambda_{2}=\psi(x) \label{s}
\end{equation} 
under the assumption that $\psi$ has bounded second derivatives. Here $\lambda_i$'s are the eigenvalues of the Hessian matrix $D^2u$ and then the phase $\psi$  becomes a potential for the mean curvature of the Lagrangian submanifold  $L=(x,Du(x))\subset \mathbb{C}^2$. 
 In two dimensions, the potential equation (\ref{s}) takes the equivalent form
\begin{equation}
    \cos\psi\Delta u +\sin\psi(\det D^2u-1)=0. \label{ss}
\end{equation}
The induced Riemannian metric $g$ can be written as 
$g=I_2+(D^2u)^2.   
$ 
In \cite[(2.19)]{HL}, the mean curvature vector $\vec{H}$ of the Lagrangian submanifold $(x,Du(x))\subset \mathbb{C}^n$ was shown to be
\begin{equation}
\vec{H}=J\nabla_g\psi \label{mean}
\end{equation}
 where $\nabla_g$ is the gradient operator for the metric $g$ and $J$ is the complex structure or the $\frac{\pi}{2}$ rotation matrix in $\mathbb{C}^n.$ Note that by our assumption on $\psi$, $|H|$ is bounded.

When the phase is constant, denoted by $c$, $u$ solves the special Lagrangian equation 
 \begin{equation}
\arctan \lambda_{1}+\arctan\lambda_2=c. \label{s1}
\end{equation}
 Equation (\ref{s1}) originates in the special Lagrangian geometry of Harvey-Lawson \cite{HL}. The Lagrangian graph $(x,Du(x))$ is special if and only if $(x,Du(x))$ is a
(volume minimizing) minimal surface in $(\mathbb{C}^n,dx^2+dy^2)$ \cite{HL}. 

The concavity of the arctangent operator plays a dominant role 
in the regularity of (\ref{s}) and (\ref{s1}). In any dimensions, if in \eqref{s1} we have critical phase $|c|=(n-2)\pi/2$ or supercritical phase $|c|>(n-2)\pi/2$, then $F(D^2u)$ has convex level sets, but it was shown by Yuan \cite{YY0} that this fails for subcritical phases $|c|<(n-2)\pi/2$.

We use $B_r$ to denote a ball of radius $r$ centered at the origin in $\mathbb R^2$ unless specified otherwise. Our main result is the following.
\begin{theorem}\label{main1}
Let $u$ be a $C^4$ solution of (\ref{s}) on $B_{R}\subset\mathbb{R}^2$ where $\psi\in C^{1,1}(B_{R})$. Then the following hold: 
\begin{align}
    |D^2u(0)|\leq C \exp [C\max_{B_R}\frac{|Du|}{R}] \text{   when } \delta\leq |\psi|\leq \frac{3\pi}{4} \label{H1}\\ 
  |D^2u(0)|\leq C \exp [C\max_{B_R}\frac{|Du|^{2}}{R^{2}}] \text{   when } |\psi|> \frac{3\pi}{4}\label{H2}
\end{align}

where $C>0$ depends on $||\psi||_{C^{1,1}(B_{R})}$, $\delta>0$, and the dimension $n=2$.
\end{theorem}
\begin{remark}\label{lem}
Note that in dimension two, the Hessian estimates derived in \cite[Theorem 1.1]{AB} hold good for all ranges of $\psi$ (including $|\psi|<\delta$) provided $\psi\in C^{1,1}(B_R)$, as discussed later in Remark (\ref{two}) of this paper. However, the existence of interior gradient estimates for solutions of (\ref{s}) when $|\psi|<\delta$ is an open problem.
\end{remark}

 In \cite{AB}, Hessian estimates for solutions of (\ref{s}) were derived in dimensions $n\geq 2$ for supercritcial $C^{1,1}$ phase. In this note, we illustrate a simplified approach in the spirit of Warren-Yuan \cite{WY2d} to prove a modified Hessian bound for solutions of (\ref{s}) in two dimensions. We provide a new
method of proof if the phase satisfies additional conditions on its size: under the assumption that $\delta\leq |\psi|\leq \frac{3\pi}{4}$ on $B_R$ we improve
the estimate in \cite{AB} to linear dependence on the gradient of the potential compared to quadratic dependence in \cite[Theorem 1.1]{AB}; under the assumption that $|\psi|>\frac{3\pi}{4}$ on $B_R$, we obtain a similar
quadratic dependence on the gradient of the potential as in \cite[Theorem 1.1]{AB} but using a simplified approach.

We observe that the previous approach in \cite{AB}, required the Michael-Simon mean value and Sobolev inequalities on generalized submanifolds of $\mathbb{R}^n$ \cite{MS}. However, in this proof,  we rely on  Warren-Yuan's \cite{WY2d} super-isoperimetric inequality, which is available only on the level sets of subharmonic functions in two dimensions. With the aid of the Jacobi inequality derived in \cite{AB}, we construct a subharmonic slope function to which we apply the super-isoperimetric inequality. This allows us to bound the  Hessian of $u$ by its integral and an integral of its gradient, and then by the volume of the Lagrangian graph. Finally we take advantage of the special two dimensional form of the equation, given by (\ref{ss}), to bound the volume element. This leads to a modified Hessian bound for (\ref{s}).

For the two dimensional case, Heinz \cite{H} derived a Hessian bound for solutions of the Monge-Ampère type equation including (\ref{s1}); Pogorelov \cite{P1} derived Hessian estimates for solutions of these equations including (\ref{s1}) with $|\psi|\geq\frac{\pi}{2}$. Gregori \cite{Gg} extended Heinz’s estimate to a gradient
bound in terms of the heights of the two dimensional minimal surfaces, and for graphs with non-zero mean curvature an additional requirement on the length of the mean curvature vector was assumed.

Higher dimensions: for critical and supercritical phases, Hessian estimates for \eqref{s1} have been obtained in \cite{WY, WdY,L}. For subcritical phases, $C^{1,\alpha}$ solutions of \eqref{s1} were constructed in \cite{NV, WaY}. For convex solutions of (\ref{s1}), a priori estimates and interior regularity were obtained in \cite{WYJ, CSY}. Recently, regularity for viscosity solutions of (\ref{s}) was studied in \cite{BS1,BS} under certain assumptions on the regularity of the phase and convexity properties of the solution. \\

 \textbf{Acknowledgments.} The author is very grateful to Yu Yuan for discussions.  The author thanks the anonymous referee for providing insightful comments and suggestions, which improved the paper.

 \section{Preliminary Inequalities}
 We will use the following results to prove higher regularity in the next section. We state the results here for the convenience of the reader.

\subsection*{Super isoperimetric inequality} 

\begin{proposition}\cite[Proposition 2.1]{WY2d}\label{supiso}
 Let $f$ be a smooth, non-negative function on $B_2\subset \mathbb{R}^2$. Suppose that $f$ satisfies the weak maximum principle: $f$ attains its maximum on the boundary of any subdomain of $B_2$. Then 
 \begin{equation}
     ||f||_{L^{\infty}(B_1)}\leq \int_{B_2}|Df|dx+\int_{B_2}f dx. \label{iso}
 \end{equation}
 \end{proposition}
 
\subsection*{Jacobi inequality}
We state the following two dimensional version of the Jacobi inequality.

\begin{lemma}\label{j1}
Let $u$ be a smooth solution of (\ref{s}) in $\mathbb{R}^{2}$ where $\psi\in C^{1,1}(B_R)$. Suppose that the Hessian $D^{2}u$ is diagonalized at $x_0$ and that the ordered eigenvalues $\lambda_1\geq\lambda_2$ of the Hessian satisfy $\lambda_1>\lambda_{2}$ at $x_0$. Then the function $b=\ln\sqrt{1+\lambda_1^2}$ is smooth near $x_0$ and at $x_0$ it satisfies 
\begin{equation}\Delta_g b\geq c(2)|\nabla_gb|^2-C\label{J}
\end{equation}
where $C=C(||\psi||_{C^{1,1}(B_1)},2)$.

\end{lemma}

\begin{proof}
The proof follows verbatim from \cite[Lemma 4.1]{AB}.
\end{proof}

\begin{remark} \label{two}
It is worth noting that the above version of the Jacobi inequality can be generalized to critical phases in any dimension $n\geq2$, i.e., the Jacobi inequality derived in \cite[Lemma 4.1]{AB} holds good when $|\psi|\geq (n-2)\frac{\pi}{2}$ and $\psi\in C^{1,1}(B_R)$ with the constant $C$ depending only on $n$ and $||\psi||_{C^{1,1}(B_R)}$. The proof follows verbatim from \cite[Lemma 4.1]{AB}.
This in turn proves that in any dimension $n\geq 2$, the Hessian estimates derived in \cite[(1.3)]{AB} hold good when $|\psi|\geq (n-2)\frac{\pi}{2}$ and $\psi\in C^{1,1}(B_R).$ However, as pointed out in Remark \ref{lem}, the existence of interior gradient estimates when $|\psi|\geq (n-2)\frac{\pi}{2}$ remains an open problem. 
\end{remark}

\section{Proof of the main result}

\begin{proof}
\noindent
Step 1. Our goal is to apply the super isoperimetric inequality to a suitable slope function. Since the slope function in Lemma \ref{j1} does not necessarily satisfy the weak maximum principle, we modify it by adding a suitable quadratic function to it. 
We assume $\psi\geq \delta$ since by symmetry $\psi\leq -\delta$ can be treated similarly.

\begin{claim}
Let $u$ be a $C^4$ solution of (\ref{s}) on $B_{R}\subset\mathbb{R}^2$ where $\psi\in C^{1,1}(B_{R})$ and $\psi\geq \delta$. Then there exists a constant $A>0$ such that the function $\Tilde{b}(x)=b(x)+\frac{A}{2}|x|^2$ is subharmonic.
\end{claim}
\begin{proof}
We assume that $D^2u$ is diagonalized at a point $p$ and $\lambda_1(p)\geq \lambda_2(p)$. Let $\lambda_1(p)>\lambda_2(p)$. Noting 
\begin{equation*}
  g^{22}=\frac{1}{1+\lambda_2^{2}}>C(\delta)>0,
\end{equation*}
and using (\ref{J}) we get 
\begin{align*}
\Delta_g \tilde{b}=\Delta_g b+\Delta_g(\frac{A}{2}|x|^2)
\geq C(2)|\nabla_gb|^2-C -\lambda_ig^{ii}\psi_i\partial_i(\frac{A}{2}|x|^2)\\
\geq -C-A\lambda_ig^{ii}\psi_i|x|]\\
\geq -C-AC_0|x|]
\end{align*}
where $0<C_0=C_0(|D\psi|)$. By scaling we can choose $|x|\leq \frac{C(\delta)}{2C_0}$ and $A=\frac{4C}{C(\delta)}$, which gives us
\begin{equation}
    \Delta_g \tilde{b}\geq 0 \label{sub}.
\end{equation}
The case $\lambda_1(p)=\lambda_2(p)$ follows similarly using \cite[Lemma 4.1]{AB}.
\end{proof}

\begin{remark}
Note that we can perform the above scaling since the constant $C_0(|D\psi|)$ rescales accordingly making $x.D\psi$ scaling invariant.
\end{remark}
\begin{remark}
Observe that the above claim holds good in dimensions $n\geq 2$ for $|\psi|\geq (n-2)\frac{\pi}{2}+\delta$. 
\end{remark}

\noindent
Step 2. We now use the special properties of the two dimensional volume element to prove the desired estimate.
For simplifying notation, we assume $R=4$ and $u$ is a solution on $B_{4}$. Then by scaling $v(x)=\frac{u(\frac{R}{4}x)}{(\frac{R}{4})^2}$, we get the estimate in Theorem \ref{main1}. 

From (\ref{sub}) we see that $\tilde{b}$ is subharmonic with respect to the induced metric on $B_2$ and hence satisfies the weak maximum principle.
So, by Lemma \ref{supiso} we get 
\begin{align*}
    ||\tilde{b}||_{L^{\infty}(B_1)}\leq \int_{B_2} |D\tilde{b}|dx+\int_{B_2}\tilde{b}dx.
\end{align*}
Denoting $\partial_i b=b_i$ we get
\begin{align}
||b||_{L^{\infty}(B_1)}\leq \int_{B_2}|Db|dx+\int_{B_2} bdx+4A|B_2|\nonumber\\
\leq \bigg(\frac{b_1^2}{1+\lambda_1^2}+\frac{b_2^2}{1+\lambda_2^2}\bigg)^{\frac{1}{2}}\bigg((1+\lambda_1^2)(1+\lambda_2^2)\bigg)^{\frac{1}{2}}dx+\int_{B_2}\frac{b}{\sqrt{\det g}}\sqrt{\det g}dx+4A|B_2|\nonumber\\
\leq [\int_{B_2} |\nabla_gb|^2dv_g]^{\frac{1}{2}}[\int_{B_2}dv_g]^{\frac{1}{2}}+\int_{B_2}dv_g+4A|B_2| \label{lala}
\end{align}
where we used $\frac{b}{\sqrt{\det g}}\leq 1$.
We choose a cut off function $\phi\in C^\infty_0(B_3)$ such that $\phi\geq 0$, $\phi=1$ on $B_2$ and $|D\phi|<2$. Using the integral version of the Jacobi inequality \cite[Proposition 4.2]{AB} we get
\begin{align}
  \int_{B_3}\phi^2|\nabla_gb|^2dv_g\leq \frac{1}{C(2)}[\int_{B_3}\phi^2\Delta_gbdv_g+\int_{B_3}\phi^2C   dv_g] \nonumber \\
    =-\frac{1}{C(2)}[\int_{B_3}\langle2\phi\nabla_g\phi,\nabla_gb\rangle dv_g+\int_{B_3}\phi^2Cdv_g] \nonumber\\
    \leq \frac{1}{2}\int_{B_3}\phi^2|\nabla_gb|^2dv_g+\frac{2}{C(2)^2}\int_{B_3}|\nabla_g\phi|^2dv_g+\frac{1}{C(2)}\int_{B_3}\phi^2Cdv_g \nonumber\\
    \implies \int_{B_2}|\nabla_gb|^2dv_g\leq
    \frac{4}{C(2)^2}\int_{B_3}|\nabla_g\phi|^2dv_g+\frac{2}{C(2)}\int_{B_3}\phi^2Cdv_g. \label{ii}
    \end{align}

     From
     \[(1+i\lambda_1)(1+i\lambda_2)=1-\sigma_2+i\sigma_1=Ve^{i\psi}
     \] we see
     $\sigma_1>0$ when $0<\psi<\pi$, and $\sigma_2>1$ when $\frac{\pi}{2}<\psi<\pi$. We use this to bound the volume element by considering the following two cases.\\
     
     \noindent
     \textbf{Case 1:} $\delta\leq\psi\leq\frac{3\pi}{4}$. \\
     We follow the argument used in \cite{WY2d} to bound the volume element.
     From (\ref{s}) we get the following simplied form of the volume element
     \[V=\bigg((1+\lambda_1^2)(1+\lambda_2^2)\bigg)^{\frac{1}{2}}=\frac{\sigma_1}{\sin\psi}.
     \]
     Since $\delta<<\frac{3\pi}{4}$, we get 
     \[V\leq \frac{1}{\sin\delta}\Delta u
     \]
and
\[\int_{B_2}dv_g\leq \frac{C(2)}{\sin\delta}||Du||_{L^{\infty}(B_3)}.
\]
Using (\ref{ss}) we observe
\begin{align*}
    |\nabla_g \phi|^2V\leq |D\phi|^2\bigg(\frac{1}{1+\lambda_1^2}+\frac{1}{1+\lambda_2^2}  \bigg)V=|D\phi|^2\bigg( \frac{2+\lambda_1^2+\lambda_2^2}{V}\bigg)\\
    =|D\phi|^2[2(1-\sigma_2)+\sigma_1^2]\frac{\sin\psi}{\sigma_1}=|D\phi|^2(2\cos\psi+\sigma_1\sin\psi).
\end{align*}
So, (\ref{ii}) reduces to \begin{align*}
    \int_{B_2}|\nabla_gb|^2dv_g\leq C[\int_{B_3}(2\cos\psi+\sigma_1\sin\psi)dx+\frac{C(2)}{\sin\delta}||Du||_{L^{\infty}(B_3)}]\\
    \leq C(2,\delta,\psi)(1+||Du||_{L^{\infty}(B_3)}).
\end{align*}

Plugging the above into (\ref{lala}), and combining the constants we get
\[||b||_{L^{\infty}(B_1)}\leq C(2,\delta,\psi)\bigg([(1+||Du||_{L^{\infty}(B_3)})||Du||_{L^{\infty}(B_3)}]^{\frac{1}{2}}+ ||Du||_{L^{\infty}(B_3)}+1 \bigg).
\]
Exponentiating, we get the estimate in (\ref{H1}).\\

\noindent
\textbf{Case 2:} $\psi>\frac{3\pi}{4}.$\\
From (\ref{lala}) and (\ref{ii}), so far we have
\begin{align}
    ||b||_{L^{\infty}(B_1)}\leq 
    C\bigg([(\int_{B_3}|\nabla_g\phi|^2dv_g+\int_{B_3}dv_g)\int_{B_3}dv_g]^{\frac{1}{2}}+\int_{B_3}dv_g+1\bigg). \label{req}
\end{align}
From \cite[Page 874]{WY2d}, for $\psi>\frac{\pi}{2}$, we see
\begin{align*}
    \int_{B_3}dv_g\leq  |\sec\psi|||Du||^2_{L^{\infty}(B_4)}.
\end{align*} Therefore for $\psi>\frac{3\pi}{4}$, we get
\[\int_{B_3}dv_g\leq \sec (\frac{3\pi}{4})||Du||^2_{L^{\infty}(B_4)}.
\]

So, (\ref{ii}) reduces to \begin{align*}
    \int_{B_2}|\nabla_gb|^2dv_g
    \leq C[\int_{B_3} (2\cos\psi+\sigma_1\sin\psi)dx+\sec (\frac{3\pi}{4})||Du||^2_{L^{\infty}(B_4)}]\\
    \leq C(2,\delta,\psi)(1+||Du||_{L^{\infty}(B_3)}+||Du||^2_{L^{\infty}(B_4)}).
\end{align*}

Plugging the above in (\ref{lala}) we get
\[||b||_{L^{\infty}(B_1)}\leq C(2,\delta,\psi) \bigg((1+||Du||_{L^{\infty}(B_3)}+||Du||^2_{L^{\infty}(B_4)})^{\frac{1}{2}}||Du||_{L^{\infty}(B_4)}+||Du||^2_{L^{\infty}(B_4)}+1\bigg).
\]

Exponentiating, we get the estimate in (\ref{H2}). 
This completes the proof of Theorem \ref{main1}. 
\end{proof}

\begin{remark}
Alternatively, one can use the Michael-Simon mean value inequality \cite[Theorem 3.4]{MS} and the $W^{1,1}$ version of the general Sobolev inequality \cite[Theorem 2.1]{MS} to derive a Hessian estimate (in terms of the gradient and the $C^{1,1}$ norm of $\psi$) for solutions of (\ref{s})  for all ranges of $\psi$ and $\psi\in C^{1,1}(B_R)$. We briefly describe the proof here: We apply a certain version of the MVI \cite[Proposition 5.1]{AB} to $b$ defined in Lemma (\ref{j1}). Next, we choose a cut off $\phi\in C^\infty_0(B_2)$ such that $\phi\geq 0$, $\phi=1$ on $B_1$ and $|D\phi|<2$. We get
\begin{equation*}b(0)\leq C\int_{\Tilde{B_1}\cap X}bdv_g\leq C\int_{B_2} b\phi dv_g
\label{bbb}
\end{equation*}
where $X=(x,Du(x))\subset \mathbb{R}^2\times\mathbb{R}^2$ is the Lagrangian submanifold, $\Tilde{B_1}$ is the ball with radius $1$ and center at $(0,Du(0))$ in $\mathbb R^2\times \mathbb R^2$, and $B_1$ is
the ball with radius $1$ and center at $0$ in $\mathbb R^2.$ 
Applying the general Sobolev inequality \cite[Theorem 2.1]{MS} on this Lagrangian submanifold, and using the mean curvature formula (\ref{mean}), we get
\begin{align*}
    b(0)\leq C[\int_{B_2}|\nabla_g(b\phi)|dv_g+\int_{B_2}|b\phi\nabla_g\psi|dv_g]\\
    \leq C\bigg[ \bigg(\int_{B_2}|\nabla_g(b\phi)|^2dv_g\bigg)^{\frac{1}{2}}\bigg(\int_{B_2}dv_g\bigg)^{\frac{1}{2}}+\bigg(\int_{B_2}|b\phi\nabla_g\psi|^2dv_g\bigg)^{\frac{1}{2}}\bigg(\int_{B_2}dv_g\bigg)^{\frac{1}{2}}\bigg].
\end{align*}
Using
\begin{align*}
    |\nabla_g(b\phi )|^2
    \leq \phi^2|\nabla_g b|^2+b^2|\nabla_g \phi|^2
    \end{align*}
the remaining proof follows from \cite{AB}.

\end{remark}

\bibliographystyle{amsalpha}
\bibliography{PJM}

\end{document}